\documentclass[12pt,twoside]{amsart}
\usepackage{amssymb}

\nonstopmode

\numberwithin{equation}{section}
\font\tengothic=eufm10 scaled\magstep 1
\font\sevengothic=eufm7 scaled\magstep 1
\newfam\gothicfam
\textfont\gothicfam=\tengothic
\scriptfont\gothicfam=\sevengothic

\newcommand {\PP}{\mathbb{P}}
\DeclareMathOperator{\pnt}{\raise 0.5mm \hbox{\large\bf.}}

\newtheorem{theorem}{Theorem}[section]

\newtheorem{proposition}[theorem]{Proposition}

\theoremstyle{definition}
\newtheorem{remark}[theorem]{Remark}
\newtheorem{example}[theorem]{Example}

\begin{document}

\title{The Geometry of the Weak Lefschetz Property and level sets of points}
\author{Juan C.\ Migliore}
\address{Department of Mathematics, University of Notre Dame,
Notre Dame, IN 46556, USA}
\email{Juan.C.Migliore.1@nd.edu}

\thanks{Part of the work for this  paper was done while the author 
was sponsored by the National Security Agency  under Grant Number
MDA904-03-1-0071. 
 }

\begin{abstract}
In a recent paper, F.\ Zanello showed that level Artinian algebras in 3 variables can fail to have the Weak Lefschetz Property (WLP), and can even fail to have unimodal Hilbert function.  We show that the same is true for the Artinian reduction of reduced, level sets of points in projective 3-space.  Our main goal is to begin an understanding of how the geometry of a set of points can prevent its Artinian reduction from having WLP, which in itself is a very algebraic notion.  More precisely, we produce level sets of points whose Artinian reductions have socle types 3 and 4 and arbitrary socle degree $\geq 12$ (in the worst case), but fail to have WLP.  We also produce a level set of points whose Artinian reduction fails to have unimodal Hilbert function; our example is based on Zanello's example.  Finally, we show that a level set of points can have Artinian reduction that has WLP but fails to have the Strong Lefschetz Property.  While our constructions are all based on basic double G-linkage, the implementations use very different methods.
\end{abstract}

%%%%%%%%%%%%%%%%%%%%%%%%%%%%%%%%%%%%%
\maketitle
\tableofcontents
%%%%%%%%%%%%%%%%%%%%%%%%%%%%%%%%%%%%%%%%%%%%

\section{Introduction}

Let $K$ be an infinite field, and let $R := K[x_0,x_1,\dots ,x_n]$ and $S := K[x_1,\dots ,x_n]$ (although for most of this paper we will use $n=3$).  For a subscheme $Z \subset \PP^n$, we denote by $h_Z$ the Hilbert function
\[
h_Z(t) = \dim (R/I_Z)_t = \dim R_t - \dim (I_Z)_t.
\]
An Artinian graded  algebra $A$ over $S$ is {\em level} if its minimal free resolution ends with a free module $S(-t)^r$ for some integers $t$ and $r$; equivalently, $A$ is level if its canonical module is generated in one degree.  The {\em socle type} of $A$ is just the integer $r$.  The {\em socle degree} of $A$ is the degree in which the socle occurs in $A$, which is the last degree in which $A$ is non-zero.  This is equal to $t-n$.  (Both the socle type and the socle degree are defined more generally, but we are only interested here in the level case.)
An arithmetically Cohen-Macaulay subscheme $V$ of $\PP^n$ is said to be {\em level} if one (hence every) Artinian reduction is level.

An Artinian graded algebra $A$ is said to have the {\em Weak Lefschetz Property} (henceforth abbreviated WLP) if, for a general linear form $L$, the induced map 
\[
\times L : A_t \rightarrow A_{t+1}
\]
has maximal rank, for all $t$.  $A$ has the {\em Strong Lefschetz Property} (SLP) if, for a general form $F$ of degree $d$, the induced map
\[
\times F : A_t \rightarrow A_{t+d}
\]
has maximal rank, for all $t$ and all $d$.

In the monograph \cite{GHMS}, the authors made an extensive study of level algebras, their properties and some methods to construct such algebras.  A study was also made of two important special cases: sets of reduced points whose Artinian reduction is level, and Artinian level algebras satisfying WLP.  Important preliminary results were obtained, and an attempt was made to bring some geometry into this very algebraic subject.  This monograph left open the following questions:

\begin{itemize}

\item Does every Artinian level algebra have the WLP?

\item Does every Artinian level algebra have unimodal Hilbert function?
\end{itemize}

Of course the primary construction to date for Artinian level algebras has been through the use of inverse systems, discovered a century ago by Macaulay \cite{orig-macaulay}.  This approach has led to wonderful results and examples, especially in the Gorenstein case, by many authors; these are too numerous to name here, and we refer to the book \cite{IK} of Iarrobino and Kanev for a detailed description.  

In a very recent paper, Zanello \cite{zanello} has completely answered the two questions listed above.  He gives examples of Artinian level algebras with unimodal Hilbert function but failing to have WLP, and he gives examples of Artinian level algebras with non-unimodal Hilbert function, and even shows that such algebras can have infinitely many ``valleys."  This is a very exciting development.
His main technique is the method of inverse systems, using at least one form chosen generically.  He also gives an example using monomial ideals, but the construction seems somewhat ad hoc.  

What is missing is the geometric side of the problem.  One of the deep questions, which was also addressed in \cite{GHMS}, is to determine which properties that do occur for Artinian algebras can be found for the Artinian reduction of a reduced set of points.  Certainly some interesting properties do not lift.  It is known, for instance, that in $k[x_1,\dots,x_6]$, a general choice of a 17-dimensional vector space of quadrics generates an Artinian ideal with Hilbert function $(1,6,4)$ whose minimal free resolution has no ghost (i.e.\ redundant) terms, while no set of 11 points in $\PP^6$  with this $h$-vector (whose ideal is generated in general by a 17-dimensional vector space of quadrics)  has such a resolution \cite{EP}.

Zanello's example using monomial ideals of course lifts using distractions (i.e.\ liftings of monomial ideals) to a reduced sets of points in $\PP^3$, so from his work we already have a reduced set of points with $h$-vector $(1,3,5,5)$ such that a particular Artinian reduction is level but does not have WLP.  It is not obvious a priori that the {\em general} Artinian reduction also fails to have WLP (which is an open condition).  See \cite{hart}, \cite{GGR}, \cite{MN2} for more on lifting monomial ideals.  

This paper represents an attempt to use a variety of geometric methods (e.g.\ union, hyperplane section) to generalize this, and to address the following questions.  The methods are almost disjoint from those of Zanello's, but they are  inspired by his beautiful paper.  

\begin{enumerate}
\item Can one recognize from the geometry of a set of points in $\PP^3$ that its general Artinian reduction is level but fails to have WLP?

\item What socle degrees and types can occur for the general Artinian reduction of a reduced set of points which is level but fails to have WLP?

\item Can the general Artinian reduction of a reduced set of points be level but not have unimodal Hilbert function?

\item Can a level, reduced set of points have the property that a particular Artinian reduction fails to have WLP, but the general Artinian reduction {\em does} have WLP?

\end{enumerate}

The author of the current paper has wondered if WLP, or at least unimodality of the $h$-vector, was one of the properties guaranteed for the Artinian reduction of a level set of points.  Zanello's results, especially the example coming from monomial ideals, began the answer.  This paper continues the study.  We note that while the underlying construction is very similar in all cases studied here, there are significant differences in the details of these proofs, and some approaches work well in some cases and fail in others.  This shows the flexibility of this approach, and leads one to hope that more such results can be obtained following this general method, for instance in higher projective spaces. Our main results are the following.  

\begin{itemize}
\item For any $d \geq 5$ there is a reduced set of points $Z \subset \PP^3$ such that {\em every} Artinian reduction of $Z$ is level of socle degree $d$ and socle type 4, but fails to have WLP.

\item For any odd $d \geq 7$ there is a reduced set of points $Z \subset \PP^3$ such that {\em every} Artinian reduction of $Z$ is level of socle degree $d$ and socle type 3, but fails to have WLP.

\item For any even $d \geq 12$ there is a reduced set of points $Z \subset \PP^3$ such that {\em every} Artinian reduction of $Z$ is level of socle degree $d$ and socle type 3, but fails to have WLP.

\item There is a reduced, level set of points in $\PP^3$ whose $h$-vector is non-unimodal.

\item There exist reduced sets of points $Z \subset \PP^3$ such that the general Artinian reduction has WLP but not SLP. 

\item Let $d \geq 2$ be a positive integer.  Then there exists a reduced, level set of points, $Z \subset \PP^3$, whose $h$-vector has a constant value in each of $d$ consecutive degrees, and for which the corresponding $d-1$ maps  in {\em any} Artinian reduction all fail to be surjective (and hence WLP fails in each of these degrees). 

\end{itemize}

It should be mentioned that Zanello found the first example of level type 3 Artinian algebras without WLP.  Our main contribution in this regard is a more detailed description of the possible socle degrees, giving all degrees beyond a certain point, and of course the fact that we get reduced sets of points.  Again, Zanello's use of inverse systems makes it impossible (with known methods) to use his main methods for points.  Another of our contributions is to show how this can be mimicked somewhat, nevertheless (see for instance Example \ref{zanello nonunimodality}).

At first glance it seems strange to ask that WLP should fail (let us say in the multiplication by a general linear form from degree $d$ to degree $d+1$) for the Artinian reduction, $A$, of the coordinate ring of a set of points, $Z$,  but the algebra $A$ still be level.  The easiest way for WLP to fail is that there is unexpected socle, but this is ruled out by the property of being level.  So instead, somehow the kernel must be bigger than expected for more subtle reasons.  How can we predict that this will happen, or explain its occurence?  In Propositions \ref{WLP criterion} and \ref{injectivity criterion} we give necessary and sufficient conditions for surjectivity and injectivity, respectively.  Here we wish to give an intuitive explanation.

The sets of points that we construct will be unions $Z = X_1 \cup Z_2$ of specially chosen sets of points, put together in a precise way by basic double G-linkage (see Proposition \ref{main tool}).  These are chosen to have $h$-vectors of a very particular sort.  It turns out that the Artinian reduction, $A$,  in degrees $d$ and $d+1$ acts as though these components decomposed into the direct sum of the Artinian reductions of $X_1$ and $X_2$, and then it is seen that the multiplication cannot have maximal rank.  For example, in Example \ref{Zanello's points} we have $Z = X_1 \cup X_2$ where the $h$-vectors are given by

 {\small 
\[
\begin{array}{l|ccccccccccccccccccccccccccccccccccc}
\deg & 0 & 1 & 2 & 3    \\
\hline
\Delta h_{X_1} && 1 & 2 & 1   \\
\Delta h_{X_2} & 1 & 2 & 3 & 4  \\ \cline{1-5}
\Delta h_Z & 1 & 3 & 5 & 5
\end{array}
\]
}

\noindent We see that even though WLP for $Z$ would require that a general linear form induces an isomorphism from degree 2 to degree 3, the decomposition gives that the rank is only $1+3 = 4$, so WLP fails.  This is not a completely accurate picture, and Example \ref{limitations} shows that we can set up situations where this philosophy predicts that WLP fails, whereas in fact it does occur.  The true explanation comes from Propositions \ref{WLP criterion} and \ref{injectivity criterion}.  Nevertheless, this intuitive approach gave us all the examples from which we produced our results.

In what way is the geometry of the set of points related to the failure of WLP?  The main tool that we use, Proposition \ref{WLP criterion}, says that the multiplication from degree $d$ to degree $d+1$ should fail to be surjective if and only if a generally chosen line $\lambda$ fails to impose $d+2$ independent conditions on $(I_Z)_{d+1}$.  We apply this to certain unions of points, where we can see from the geometry that $\lambda$ has to impose fewer than the expected number of conditions because of fixed 2-dimensional components in the base locus of the linear system.  Proposition \ref{injectivity criterion} gives a condition that is equivalent to the failure of injectivity, and we illustrate its use in Example \ref{not injective}.

In this paper we only work in codimension three, studying points in $\PP^3$, but the methods generalize at least in part to any $\PP^n$.  Many of our examples were verified with the computer program CoCoA \cite{cocoa}.

The author is very grateful to Fabrizio Zanello for sending an early version of his inspiring paper.  He is also grateful to Tony Geramita and to Uwe Nagel for very enjoyable conversations about Zanello's work, which led to the  ideas developed here.

\newpage

%%%%%%%%%%%%%%%%%%%%%%%%%%%%%%%%%%%%%%%%%%%%%%

\section{The machinery}

The purpose of this paper is to produce, by very geometric methods (e.g.\ unions, hyperplane and hypersurface sections) reduced sets of points whose Artinian reduction is level but does not have the Weak Lefschetz property.  In this section we discuss how our sets are constructed, give the reason why they are level, and show how the geometry of the points allows an analysis to help determine if the Weak Lefschetz property holds or not.  These ideas are applied to specific questions in the subsequent sections.

\begin{proposition}[\cite{KMMNP}] \label{main tool}
Let $C \subset \PP^3$ be an arithmetically Cohen-Macaulay curve and let $Z_1 \subset C$ be any zero-dimensional scheme.  Let $F$ be a non-zero divisor on  $R/I_{C}$ of degree $d$, and let $Z_2$ be the hypersurface section of $C$ by the surface $S$ defined by $F$.  Let $Z$ be the zero-scheme defined by the ideal $I := L \cdot I_{Z_1} + I_C$.  Then

\begin{itemize}
\item[(a)] We have the exact sequence 
\[
0 \rightarrow I_C(-d) \rightarrow I_{Z_1}(-d) \oplus I_C \rightarrow I_Z \rightarrow 0
\]

\item[(b)] $I$ is saturated, i.e.\ $I = I_Z$.

\item[(c)] We have the exact sequence
\[
0 \rightarrow (R/I_C)(-1) \rightarrow (R/I_{Z_1})(-1) \oplus (R/I_C) \rightarrow R/I_Z \rightarrow 0.
\]

\item[(d)] The Hilbert function of $Z$ is given by the formula
\[
\begin{array}{rcl}
h_Z (t) & = & h_{Z_1}(t-d) + h_C(t) - h_C(t-d) \\
& = & h_{Z_1}(t-d) + h_{Z_2}(t)
\end{array}
\]

\item[(e)] If $Z_1$ and $C$ are both level, and if they have the same shift, $-t$,in the last free modules of their minimal free resolutions (which means that the socle degrees  of their Artinian reductions differ by one since they have different codimensions), then $Z$ is level with minimal free resolution having shift $-t-d$.  The socle type is less than or equal to the sum of the socle types of $C$ and of $Z_1$.

\item[(f)] If in addition $F$ is a non-zero divisor on $R/I_{Z_1}$ (i.e.\ $Z_1$ and $Z_2$ are disjoint) then $Z = Z_1 \cup Z_2$ as schemes.
\end{itemize}
\end{proposition}

\begin{proof}
Part (a) can be found in \cite{KMMNP}, Lemma 4.8, and the rest follows immediately from this.  Part (e), in particular, follows from (a) using a mapping cone.
\end{proof}

\begin{remark} \label{common use}
\begin{enumerate}
\item Following \cite{KMMNP}, we will refer to this construction as {\em basic double G-linkage}.  
Proposition \ref{main tool} is often used in this paper in the following way.  Let $C_1 \subset C_2$ be arithmetically Cohen-Macaulay curves in $\PP^3$.  Then the general hypersurface section of $C_1$ of any degree $d_1$ is also a subscheme of $C_2$, and so the union of the general degree $d_1$ hypersurface section of $C_1$ and the general degree $d_2$ hypersurface section of $C_2$ fits into the scheme of Proposition \ref{main tool}.  These curves, and the degrees of the hypersurfaces, can be chosen so that the union is level, and then geometric properties can be read off.

\item In this paper we will use the first difference of the Hilbert function formula above, to compute the $h$-vectors.
\end{enumerate}
\end{remark}

\begin{proposition} \label{WLP criterion}
Let $Z \subset \PP^3$ be a zero-dimensional scheme.  Let $A$ be the Artinian reduction of $R/I_Z$ by a linear form $L_1$ and let $L_2$ be a second, {\em general} linear form.  Then the multiplication map $\times L_2 : A_{t-1} \rightarrow A_t$ is surjective if and only if the line $\lambda$ defined by $(L_1,L_2)$ imposes $t+1$ conditions on the linear system $|(I_Z)_t|$.
\end{proposition}

\begin{proof}
Note that $L_1$ only has to avoid the points of $Z$.  No additional generality is required.  We have $A = R/(I_Z + (L_1))$, and the exact sequence
\[
[R/(I_Z + (L_1))]_{t-1} \stackrel{\times L_2}{\longrightarrow} [R/(I_Z + L_1))]_t \rightarrow
[R/(I_Z + I_\lambda)]_t \rightarrow 0.
\]
Hence the multiplication $\times L_2$ is surjective if and only if $(I_Z + I_\lambda)_t = R_t$.  Now we also have the exact sequence
\[
0 \rightarrow (I_Z \cap I_\lambda)_t \rightarrow (I_Z)_t \oplus (I_\lambda)_t \rightarrow (I_C + I_\lambda)_t \rightarrow 0.
\]
The statement that $\lambda$ imposes $t+1$ conditions on $|(I_Z)_t|$ is equivalent to the statement that $\dim (I_Z \cap I_\lambda)_t = \dim (I_Z)_t - (t+1)$.   Since $\dim (I_\lambda)_t = \dim R_t - (t+1)$, 
the statement now follows immediately from exactness of the sequence.
\end{proof}

\begin{example}\label{Zanello's points}
We illustrate the methods that will be used in this paper by constructing a level reduced set of points with $h$-vector $(1,3,5,5)$ whose Artinian reduction (arbitrary) does not have WLP.  The existence of such a set of points (almost) follows also from Zanello's work \cite{zanello}: his Example 7 gives an Artinian monomial ideal with this Hilbert function, and as mentioned above we can lift this example to a reduced set of points.  The only thing missing is the question of whether WLP fails for all Artinian reductions or only a special one.  In any case, it was an attempt to understand this example from a geometric point of view that motivated our current work.

Let $C$ be a smooth arithmetically Cohen-Macaulay curve with $h$-vector $(1,2,3,4)$.  A general hyperplane section of $C$ contains four points, $X_1$, in linear general position.  Note that $X_1$ is a complete intersection of type $(1,2,2)$ in $\PP^3$.  Let $X_2$ be another general hyperplane section of $C$, and let $Z = X_1 \cup X_2$.  $Z$ is clearly reduced.  By Proposition \ref{main tool}, the $h$-vector of $Z$ is computed by
 {\small 
\[
\begin{array}{l|ccccccccccccccccccccccccccccccccccc}
\deg & 0 & 1 & 2 & 3    \\
\hline
\Delta h_{X_1} && 1 & 2 & 1   \\
\Delta h_{X_2} & 1 & 2 & 3 & 4  \\ \cline{1-5}
\Delta h_Z & 1 & 3 & 5 & 5
\end{array}
\]
}
The fact that $Z$ is level is computed from an easy mapping cone using the short exact sequence
\[
0 \rightarrow I_C(-1) \rightarrow I_{X_1}(-1) \oplus I_C \rightarrow I_Z  \rightarrow 0,
\]
which gives the resolution
\[
0 \rightarrow
R(-6)^5 \rightarrow 
\left (
\begin{array}{c}
R(-5)^{10} \\
\oplus \\
R(-4)^2
\end{array}
\right )
\rightarrow 
\left (
\begin{array}{c}
R(-4)^5 \\
\oplus \\
R(-3)^2 \\
\oplus \\
R(-2)
\end{array}
\right )
\rightarrow I_Z \rightarrow 0.
\]
The failure of WLP to hold comes from Proposition \ref{WLP criterion}.  Indeed, note that any polynomial of degree 3 containing $Z$ has to contain the plane of $X_2$ as a factor, since the points of $X_2$ do not lie on any cubic on that plane.  Now let $\lambda$ be a general line, meeting the plane of $X_2$ in a point, $P$.  The point $P$  is thus a base point in the linear system $|(I_Z)_3|$, so $\lambda$ imposes at most three conditions on cubics, and hence WLP fails.

We remark that the resolution above is also the minimal free resolution of a set of 14 general points on a smooth arithmetically Cohen-Macaulay curve of degree 5 and genus 2, the Artinian reduction of which is also level, and {\em does} have WLP.    \qed
\end{example}

In this paper we are mostly concerned with the failure of WLP occurring by having the multiplication by a general linear form fail to be surjective.  However, of course it can also be injectivity that fails, and there is a criterion for that as well.  Although we will not use this criterion for any of our main results, we record it here.  

\begin{proposition} \label{injectivity criterion}
Let $Z \subset \PP^3$ be a zero-dimensional scheme.  Let $A$ be the Artinian reduction of $R/I_Z$ by a linear form $L_1$ and let $L_2$ be a second, {\em general} linear form.  Let $\lambda$ be the line defined by the ideal $(L_1,L_2)$. Then the multiplication map $\times L_2 : A_{t-1} \rightarrow A_t$ is injective if and only if $(I_Z \cdot I_\lambda)_t = (I_Z \cap I_\lambda)_t$.
\end{proposition}

\begin{proof}
Note that again, $L_1$ only has to avoid the points of $Z$.  Let 
\[
J = \frac{I_Z}{(L_1 \cdot I_Z)} \cong \frac{I_Z + (L_1)}{(L_1)}.
\]
  If $S := R/(L_1)$ then we have $A = S/J$.  By \cite{migbook} Proposition 2.1.5 (see also \cite{submod}), 
\[
\frac{J : L_2}{J}(-1) \cong \frac{I_Z \cap I_\lambda}{I_Z \cdot I_\lambda}.
\]
The result then follows from the exact sequence
\[
0 \rightarrow \frac{J : L_2}{J}(-1) \rightarrow S/J(-1) \stackrel{\times L_2}{\longrightarrow} S/J \rightarrow R/(I_Z + I_\lambda) \rightarrow 0.
\]
\end{proof}

We will illustrate the use of Proposition \ref{injectivity criterion}  in Section \ref{other unusual beh} (see Example \ref{not injective}), but we have not found an effective way to apply it for broader results.

%%%%%%%%%%%%%%%%%%%%%%%%%%%%%%%%%%%%%%%%%%%%%%%%

\section{Small socle type}

In this section we show how our methods give classes of reduced sets of points with small socle type, and how the geometry of these sets is what determines the fact that the WLP fails for the Artinian reduction.

\begin{theorem} \label{type 4}
For any $d \geq 5$ there is a reduced set of points $Z \subset \PP^3$ such that {\em any} Artinian reduction of $Z$ is level of socle degree $d$ and socle type 4, but fails to have the Weak Lefschetz Property.
\end{theorem}

\begin{proof}
We begin by showing the existence of certain level sets of points, $Z_d$, in $\PP^2$ of socle type 3 and socle degree $d$ that contain as subsets certain complete intersection sets of points $Y_d$ of socle degree $d-1$.  These can be obtained as liftings (also known as {\em distractions}) of monomial ideals in two variables, so it is enough to present these monomial ideals.   We distinguish the cases $d$ even and $d$ odd.  We first give the construction and collect some facts for both of these cases.

\bigskip

\noindent {\bf Case 1:} $d$ odd, $d = 2r-1$. \hskip 1in {\bf Case 2:} $d$ even, $d = 2r$.

%\begin{figure}[ht]
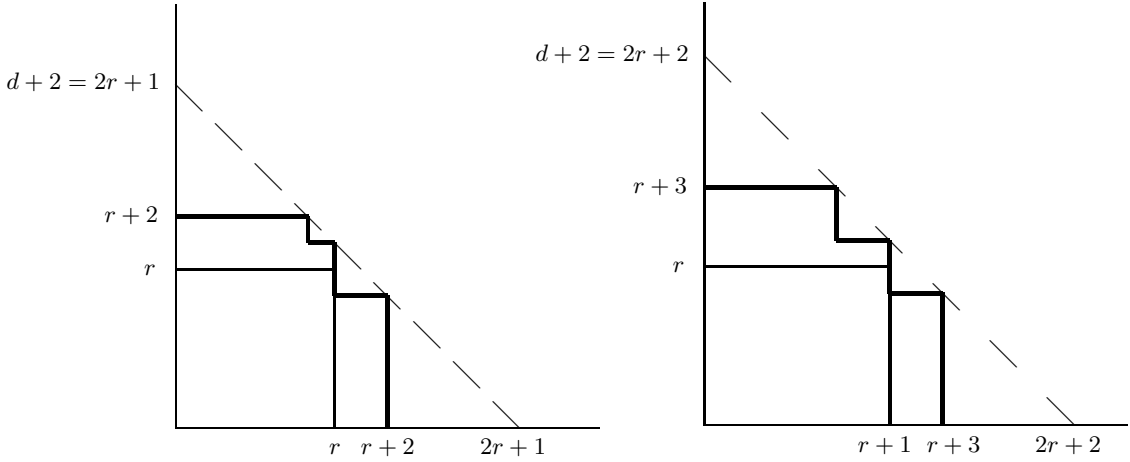
\begin{picture}(500,200)
\linethickness{0.15mm}
\put (70,30){\line(1,0){160}}
\put (70,30){\line(0,1){160}}
\linethickness{0.3mm}
\put(130,30){\line(0,1){60}}
\put(70,90){\line(1,0){60}}
\linethickness{0.5mm}
\put(70,110){\line(1,0){50}}
\put(120,110){\line(0,-1){10}}
\put(120,100){\line(1,0){10}}
\put(130,100){\line(0,-1){20}}
\put(130,80){\line(1,0){20}}
\put(150,80){\line(0,-1){50}}
\linethickness{0.15mm}
\multiput(70,160)(15,-15){9}{\line(1,-1){10}}
\put(58,88){\scriptsize $r$}
\put(43,108){\scriptsize $r+2$}
\put(6,158){\scriptsize $d+2 = 2r+1$}
\put(128,20){\scriptsize $r$}
\put(140,20){\scriptsize $r+2$}
\put(185,20){\scriptsize $2r+1$}
\end{picture}
%\caption{Monomial Ideal $M_d$ ($d$ odd) \label{odd d}}
%\end{figure}

%\begin{figure}[ht]
\begin{picture}(500,0)
\linethickness{0.15mm}
\put (270,45){\line(1,0){160}}
\put (270,45){\line(0,1){160}}
\linethickness{0.3mm}
\put(340,45){\line(0,1){60}}
\put(270,105){\line(1,0){70}}
\linethickness{0.5mm}
\put(270,135){\line(1,0){50}}
\put(320,135){\line(0,-1){20}}
\put(320,115){\line(1,0){20}}
\put(340,115){\line(0,-1){20}}
\put(340,95){\line(1,0){20}}
\put(360,95){\line(0,-1){50}}
\linethickness{0.15mm}
\multiput(270,185)(21.6,-21.6){7}{\line(1,-1){10}}
\put(258,103){\scriptsize $r$}
\put(243,133){\scriptsize $r+3$}
\put(206,183){\scriptsize $d+2 = 2r+2$}
\put(328,35){\scriptsize $r+1$}
\put(354,35){\scriptsize $r+3$}
\put(395,35){\scriptsize $2r+2$}
\end{picture}
%\caption{Monomial Ideal $M_d$ ($d$ even) \label{even d}}
%\end{figure}

In Case 1, $Z_d$ is the reduced set of points in $\PP^2$ obtained by lifting the monomial ideal (in two variables) $(y^{r+2}, x^{r-1}y^{r+1}, x^ry^{r-1}, x^{r+2})$.   This ideal is illustrated above by the thickest lines.  The complete intersection $Y_d \subset Z_d$ is obtained by lifting the complete intersection $(x^r,y^r)$, illustrated above by the second thickest lines.  $Y_d$ is a complete intersection of type $(r,r)$.  
In Case 2, $Z_d$ is obtained by lifting the monomial ideal $(y^{r+3}, x^{r-1}y^{r+1}, x^{r+1}y^{r-1}, x^{r+3})$, and $Y_d \subset Z_d$ is obtained by lifting the complete intersection $(x^{r+1}, y^r)$.  $Y_d$ is a complete intersection of type $(r,r+1)$.

The following facts are easy to verify.  For  $r=3$ some of the higher values of the Hilbert function are not correct below, but in the critical range where we have to show the failure of WLP, they are correct.

\noindent {\bf Case 1:}

\begin{itemize}
\item[(a)] $Y_d$ has $h$-vector
 {
\[
\begin{array}{l|ccccccccccccccccccccccccccccccccccc}
\deg & 0 & 1 & 2  & \dots & r-2 & r-1 & r & \dots & 2r-4 & 2r-3 & 2r-2  \\
\hline
\Delta h_{Y_d} & 1 & 2 & 3 & \dots & r-1 & r & r-1 & \dots & 3 & 2 & 1
\end{array}
\]
}
In particular, the socle degree is $2r-2 = d-1$.

\item[(b)] $Z_d$ is level with socle type 3 and socle degree $2r-1 = d$.  It has minimal free resolution (over $S = K[x_1,x_2,x_3]$) which can be read immediately from the diagram above:
\[
0 \rightarrow S(-2r-1)^3 \rightarrow 
\left (
\begin{array}{c}
S(-r-2)^2 \\
\oplus \\
S(-2r+1)  \\
\oplus \\
S(-2r)
\end{array}
\right )
\rightarrow I_{Z_d} \rightarrow 0.
\]

\item[(c)] $Z_d$ has $h$-vector
 {
\[
\begin{array}{l|ccccccccccccccccccccccccccccccccccc}
\deg & 0 & 1 & 2  & \dots & r+1 & r+2 & r+3 & \dots & 2r-3 & 2r-2 & 2r-1  \\
\hline
\Delta h_{Z_d} & 1 & 2 & 3 & \dots & r+2 & r +1& r & \dots & 6 & 5 & 3
\end{array}
\]
}
\end{itemize}

\noindent {\bf Case 2:}

\begin{itemize}
\item[(a)] $Y_d$ has $h$-vector
 {\
\[
\begin{array}{l|ccccccccccccccccccccccccccccccccccc}
\deg & 0 & 1 & 2  & \dots & r-2 & r-1 & r & r+1 & \dots & 2r-3 & 2r-2 & 2r-1  \\
\hline
\Delta h_{Y_d} & 1 & 2 & 3 & \dots & r-1 & r & r & r-1& \dots & 3 & 2 & 1
\end{array}
\]
}
In particular, the socle degree is $2r-1 = d-1$.

\item[(b)] $Z_d$ is level with socle type 3 and socle degree $2r = d$.  It has minimal free resolution (over $S = K[x_1,x_2,x_3]$) which can be read immediately from the diagram above:
\[
0 \rightarrow S(-2r-2)^3 \rightarrow 
\left (
\begin{array}{c}
S(-r-3)^2 \\
\oplus \\
S(-2r)^2
\end{array}
\right )
\rightarrow I_{Z_d} \rightarrow 0.
\]

\item[(c)] $Z_d$ has $h$-vector
 {
\[
\begin{array}{l|ccccccccccccccccccccccccccccccccccc}
\deg & 0 & 1 & 2  & \dots & r+2 & r+3 & r+4 & \dots & 2r-2 & 2r-1 & 2r  \\
\hline
\Delta h_{Z_d} & 1 & 2 & 3 & \dots & r+3 & r +2& r+1 & \dots & 7 & 6 & 3
\end{array}
\]
}
\end{itemize}

Choose a general point in $\PP^3$ and consider the cones $B_d \subset C_d$ in $\PP^3$ over $Y_d \subset Z_d \subset \PP^2$.  $B_d$ and $C_d$ are arithmetically Cohen-Macaulay curves, and have the same $h$-vector and graded Betti numbers (over $R = k[x_0,x_1,x_2,x_3]$) as $Y_d$ and $Z_d$ do, respectively, over $S$.

Let $Z$ be the union of a general hyperplane section, $X_1$, of $B_d$ and a general hyperplane section, $X_2$, of $C_d$ (see Remark \ref{common use}).  We will show that $Z$ satisfies all the claims in the statement of the theorem.  {\em We will prove it for $d$ odd ({\bf Case 1}), and leave it to the reader to check the case $d$ even, which is completely analogous.}

We first compute the $h$-vector. Observe that $X_1$ (resp.\ $X_2$) has the same $h$-vector as $Y_d$ (resp.\ $Z_d$), and that in particular $X_1 \subset C_d$.  Hence $Z$ is obtained by basic double G-linkage, and if $L$ is the linear form defining the hyperplane of $X_2$ then we have
\[
I_Z = L \cdot I_{X_1} + I_{C_d}.
\]
The $h$-vector of $Z$ is thus obtained by the computation
 {\small 
\[
\begin{array}{l|ccccccccccccccccccccccccccccccccccc}
\deg & 0 & 1 & 2 & 3 &  \dots & r-1 & r & r+1 & r+2 & \dots  & 2r-3 & 2r-2 & 2r-1 &
2r \\
\hline
\Delta h_{X_1} & & 1 & 2 & 3 & \dots & r-1 & r & r-1& r-2 &  \dots & 3 & 2 & 1 \\
\Delta h_{X_2} & 1 & 2 & 3 & 4 & \dots & r & r+1 & r+2 & r+1 & \dots  & 6 & 5 & 3 \\ \cline{1-15}
\Delta h_Z & 1 & 3 & 5 & 7 & \dots & 2r-1 & 2r+1 & 2r+1 & 2r-1 & \dots & 9 & 7 & 4
\end{array}
\]
}
(Notice the shift in $\Delta h_{X_1}$.)  In particular, note that if $T$ is the ring of an Artinian reduction and $J$ is the corresponding ideal, then $\dim (T/J)_r = \dim (T/J)_{r+1} = 2r+1$.  This is where we will focus our attention to show the failure of WLP.

We first check that $Z$ is level, of type 4.  This is immediate from the mapping cone associated to the diagram

\vskip .8cm

 { \scriptsize 
\[
\begin{array}{ccccrclccccccccccccccccccccccc}
&&&&& 0 \\
&&&&& \downarrow \\
&& 0 && R(-2r-2) & \oplus \\
&& \downarrow &&& \downarrow \\
&& R(-2r-2)^3 && R(-r-2)^2 \oplus R(-2r-1)  & \oplus & R(-2r-1)^3 \\
&& \downarrow &&& \downarrow \\
&& R(-r-3)^2 \oplus R(-2r) \oplus R(-2r-1) && R(-2) \oplus R(-r-1)^2 & \oplus & R(-r-2)^2 \oplus R(-2r+1) \oplus R(-2r) \\
&& \downarrow &&& \downarrow \\
&& 0  \rightarrow \hbox{\hskip 1.5cm}  I_{C_d}(-1) \phantom{0  \rightarrow \hbox{\hskip 1.5cm} } & \rightarrow &  I_{X_1}(-1) & \oplus & I_{C_d} \\
&& \downarrow &&&  \downarrow \\
&& 0 &&&   0
\end{array}
\]
} 
which gives a minimal (in this case) free resolution of $I_Z$.

We now verify the failure of WLP to hold for any Artinian reduction of $Z$. By Proposition \ref{WLP criterion}, it is enough to show that a general line, $\lambda$, imposes fewer than the expected $r+2$ conditions on $|(I_Z)_{r+1}|$.

Let $I_\lambda = (L_1,L_2)$ as in Proposition \ref{WLP criterion}.  By assumption, $L_1$ avoids the points of $Z$ and $L_2$ is general.  Hence $\lambda$ meets the plane of $X_2$ at a single point, $P$, disjoint from $X_2$.  Now, thanks to the $h$-vector computation of $Z_d$ (which is the same as that of $X_2$), any form of degree $r+1$ containing $Z$ must contain the plane of $X_2$ as a factor, since otherwise it would restrict to a curve of degree $r+1$ containing $X_2$.  Hence $P$ is a base point of the linear system $|(I_Z)_{r+1}|$.  Therefore $\lambda$ imposes at most $r+1$ conditions on $|(I_Z)_{r+1}|$, and we have shown that the map induced by $L_2$ on the Artinian reduction, $T$, by $L_1$ fails to be surjective, by Proposition \ref{WLP criterion}.  Since $\dim (T/J)_r = \dim (T/J)_{r+1} = 2r+1$, this means that WLP fails for $Z$.

One checks that this works for $r \geq 3$ in both the even and odd cases, completing the proof.
\end{proof}

\begin{theorem} \label{type 3 odd}
For any {\em odd} $d \geq 7$, there is a reduced set of points $Z \subset \PP^3$ such that {\em any} Artinian reduction of $Z$ is level of socle degee $d$ and socle type 3, but fails to have the Weak Lefschetz Property.
\end{theorem}

\begin{proof}
The proof follows along very similar lines to that of Theorem \ref{type 4}, so we will highlight primarily the differences.  Let $d = 2r-1$.  We again use basic double G-linkage to produce $Z$ as the union of subsets $X_1, X_2$, both degenerate but on different planes, so that $Z$ is level of type 3, and the $h$-vector of $Z$ is computed by

 {\scriptsize
\[
\begin{array}{l|ccccccccccccccccccccccccccccccccccc}
\deg & 0 & 1 & 2 & 3 &  \dots & r-2 & r-1 & r & r+1 & r+2 & \dots  & 2r-4 & 2r-3 & 2r-2 & 2r-1 \\
\hline
\Delta h_{X_1} & & 1 & 2 & 3 & \dots & r-2 & r-1 & r& r-1 & r-2 &  \dots & 4 & 3 & 2 & 1 \\
\Delta h_{X_2} & 1 & 2 & 3 & 4 & \dots & r-1 & r & r+1 & r+2 & r+1 & \dots & 7 & 6 & 4 & 2 \\ \hline
\Delta h_Z & 1 & 3 & 5 & 7 & \dots & 2r-3 & 2r-1 & 2r+1 & 2r+1 & 2r-1 & \dots & 11 & 9 & 6 & 3
\end{array}
\]
}
\!\!Again, WLP will be shown to fail by the failure of the multiplication by a general linear form from degree $r$ to degree $r+1$ to be surjective.

To achieve this, we start by constructing a set of points $Z_d \subset \PP^2$ that is level of type 2 and contains a subset, $Y_d$, that is a complete intersection of type $(r,r)$.  However, this cannot be achieved by monomial ideals, and we instead need a completely different argument.

As a first step, we show that if $Y_0$ is a general complete intersection in $\PP^2$ of type (3,3) and $Y_d$ is a general complete intersection of type $(r,r)$, then there is a reduced complete intersection of type $(r+2,r+2)$ containing $Y_0 \cup Y_d$.  To see this, consider the exact sequence 
\[
0 \rightarrow [I_{Y_0} \cap I_{Y_d}]_{r+2} \rightarrow [I_{Y_0}]_{r+2} \oplus [I_{Y_d}]_{r+2} \rightarrow [I_{Y_0} + I_{Y_d}]_{r+2} \rightarrow 0.
\]
It is easy to see that $\dim [I_{Y_0}]_{r+2} = \binom{r+4}{2} - 9$ and that $\dim [I_{Y_d}]_{r+2} = 12$.  On the other hand, $I_{Y_0} + I_{Y_d}$ is a general choice of forms of degrees $3,3,r,r$.  By a result of Anick \cite{anick} we know the Hilbert function of an ideal of general forms of any specified degree in three variables, and  we can verify that $\dim [I_{Y_0} + I_{Y_d}]_{r+2} = \binom{r+4}{2}$ (i.e.\ the Hilbert function is zero in degree $r+2$).  This gives
\[
\dim [I_{Y_0} \cap I_{Y_d}]_{r+2} = \binom{r+4}{2} - 9 + 12 - \binom{r+4}{2} = 3.
\]

We claim that a general choice of two forms in $[I_{Y_0} \cap I_{Y_d}]_{r+2}$ will define a {\em reduced} complete intersection of type $(r+2,r+2)$.  Note that to specify a general complete intersection of type $(3,3)$, it is equivalent to specify eight general points in $\PP^2$; then the ninth point is uniquely determined.  It is clear that a general choice of two forms in $[I_{Y_d}]_{r+2}$ (which has dimension 12) defines a reduced complete intersection, and that imposing the further passing through eight general points drops the dimension to 4 and preserves the property that a general choice of two defines a reduced complete intersection.  Now the ninth point of $Y_0$ is determined, but it still imposes a new condition on  $[I_{Y_d}]_{r+2}$ by the computation above.  Since everything has been chosen generically so far, we preserve the property of reducedness.

Hence by taking two general elements of the linear system $[I_{Y_0} \cap I_{Y_d}]_{r+2}$, we can link $Y_0$ to a reduced set of points $Z_d$ that is level of type 2, has $h$-vector given by the line $\Delta H_{X_2}$ in the computation at the beginning of this proof, and contains the complete intersection $Y_d$ of type $(r,r)$.  By the general choices, everything is reduced.

We then proceed as in Theorem \ref{type 4}, taking cones over $Y_d \subset Z_d$ to produce arithmetically Cohen-Macaulay curves curves $B_d \subset C_d \subset \PP^3$ over $Y_d \subset Z_d \subset \PP^2$, and form the union $Z = X_1 \cup X_2$, where $X_1$ is a general hyperplane section of $B_d$ and $X_2$ is a general hyperplane section of $C_d$.  The verification of the Hilbert function, the level property and the WLP is identical to that in Theorem \ref{type 4} and is left to the reader.

We remark that for even socle degree this method does not yield a complete intersection containing the analogous $Y_0 \cup Y_d$, and so it fails to work.  In the next result we will use a different construction for the even case.
\end{proof}

\begin{theorem} \label{type 3 even}
For any {\em even} $d \geq 12$, there is a reduced set of points $Z \subset \PP^3$ such that {\em any} Artinian reduction of $Z$ is level of socle degree $d$ and socle type 3, but fails to have the Weak Lefschetz Property.
\end{theorem}

\begin{proof}
The methods are rather different, even if the underlying construction is based on basic double G-linkage.  We will use basic double G-linkage to produce $Z$ as the union of subsets $X_1,X_2$.  This time, however, neither is degenerate.  One, $X_1$, will be arithmetically Gorenstein (not a complete intersection), and the other, $X_2$, will be a quadric hypersurface section of a suitable arithmetically Cohen-Macaulay curve in $\PP^3$.  
Another difference is that this time it will not be enough to take a cone over a set of points in $\PP^2$.  Rather, we will apply a direct geometric construction, and from that use liaison methods to produce the reduced, arithmetically Gorenstein set of points that we need.

Let $d = 2r$.  We will have to construct an arithmetically Cohen-Macaulay curve with the following properties:

\begin{enumerate}
\item \label{stick req} $C$ is a stick figure, i.e.\ it is a union of lines with at most two passing through any given point.

\item \label{subcurve req} $C$ contains as a subset an arithmetically Cohen-Macaulay curve  $Y$ with maximal $h$-vector
\[
(1,2,3,\dots,r-1,r).
\]
We will construct on $Y$ (and hence on $C$) a reduced, arithmetically Gorenstein set of points, $X_1$, with properties described below.

\item Any Artinian reduction of $R/I_C$ is a compressed level algebra with socle type 2 and socle degree $2r-1$.  Any Artinian reduction of $R/I_{X_1}$ is a compressed Gorenstein algebra with socle degree $2r-2 = d-2$.  (But notice that $C$ has codimension two and $X_1$ has codimension three.)
\end{enumerate}

Once we have these basic objects, we will construct $Z$ with basic double G-linkage as follows: $Z = X_1 \cup X_2$, where $X_2$ is a general {\em quadric} hypersurface section of $C$.  

First let us construct $C$.  We have three cases: $r \equiv 0 \hbox{ (mod 3)}, r \equiv 1 \hbox{ (mod 3)},$ and $r \equiv 2 \hbox{ (mod 3)}$.  In all cases, $C$ will be residual to a complete intersection, $D$, of type $(t,t)$, where
\[
t = 
\left \{
\begin{array}{ll}
\frac{2}{3} r & \hbox{if } r \equiv 0 \hbox{ (mod 3)}; \\ \\
\frac{2}{3} (r-1)+1 & \hbox{if } r \equiv 1 \hbox{ (mod 3)};\\ \\
\frac{2}{3}(r-2)+1 & \hbox{if } r \equiv 2 \hbox{ (mod 3)},
\end{array}
\right.
\]
inside a complete intersection $E$ of type 
\[
\left \{
\begin{array}{ll}
 \left (\frac{4r}{3},\frac{4r}{3}+1 \right ) & \hbox{if } r \equiv 0 \hbox{ (mod 3)}; \\ \\
\left (\frac{4(r-1)}{3}+2, \frac{4(r-1)}{3}+2 \right ) & \hbox{if } r \equiv 1 \hbox{ (mod 3)}; \\ \\
\left (\frac{4(r-2)}{3}+3, \frac{4(r-2)}{3}+3 \right ) & \hbox{if } r \equiv 2 \hbox{ (mod 3)}.
\end{array}
\right.
\]

From a standard mapping cone argument it is easy to see that such a $C$ is level of type 2, since it is linked to a complete intersection of balanced type $(t,t)$.  We will need information about the Hilbert function of such a $C$.  We will work it out for $r \equiv 0 \hbox{ (mod 3)}$, and leave the remaining cases to the reader.  We now compute the $h$-vector, using the formula for Hilbert functions under linkage \cite{DGO}, \cite{migbook}.  (We read the $h$-vector of $D$ from right to left, noting that it is symmetric, and that of $C$ from left to right.)

 { 
\[
\begin{array}{l|ccccccccccccccccccccccccccccccccccc}
\deg & 0 & 1 & 2 &   \dots & \frac{4r}{3}-2 & \frac{4r}{3}-1 & \frac{4r}{3} & \frac{4r}{3}+1 & \frac{4r}{3}+2 & \dots  \\
\hline
E & 1 & 2 & 3  & \dots &  \frac{4r}{3}-1 & \frac{4r}{3} & \frac{4r}{3} & \frac{4r}{3}-1 & \frac{4r}{3}-2 & \dots  \\
D & &&&  & &  & 0 & 1 & 2 & \dots \\ \cline{1-11}
C & 1 & 2 & 3 & \dots & \frac{4r}{3}-1 & \frac{4r}{3} & \frac{4r}{3} & \frac{4r}{3}-2 & \frac{4r}{3}-4 & \dots 
\end{array} 
\hbox{\hskip 1.4cm}
\]

\bigskip

\[
\begin{array}{l|cccccccccccccccccccccccccccccccccccc}
\deg & \dots &  2r-2 & 2r-1 & 2r & 2r+1 & \dots & \frac{8r}{3}-3 & \frac{8r}{3}-2 & \frac{8r}{3}-1 \\
\hline
E && \frac{2r}{3}+2 & \frac{2r}{3}+1 & \frac{2r}{3} & \frac{2r}{3}-1& \dots & 3 & 2 & 1 \\
D && \frac{2r}{3}-2 & \frac{2r}{3}-1 & \frac{2r}{3} & \frac{2r}{3}-1 & \dots & 3 & 2 & 1 \\ \hline
C & \dots & 4 & 2
\end{array}
\]
}

Now, we want $C$ to contain as a subcurve an arithmetically Cohen-Macaulay curve with maximal $h$-vector $(1,2,3,\dots,r)$.  In order to be able to construct $C$ with the above numerical properties and with a subcurve having this $h$-vector (and satisfying the other requirements listed above), we use a trick similar to that used in \cite{MN3}, and related to liftings of monomial ideals.  That is, consider two families of planes in $\PP^3$, $L_1,\dots,L_{\frac{4r}{3}}$ and $M_1,\dots,M_{\frac{4r}{3}+1}$, where the $L_i$ and the $M_i$ are chosen generically.  In the following diagram, each horizontal line represents an $L_i$ and each vertical line represents an $M_j$, and the intersection point represents the line of intersection of $L_i$ and $M_j$.  It is clear that if we denote by $L := \prod L_i$ and $M := \prod M_j$ then $(L,M)$ is a complete intersection of type $(\frac{4r}{3}, \frac{4r}{3}+1)$ which is a stick figure.

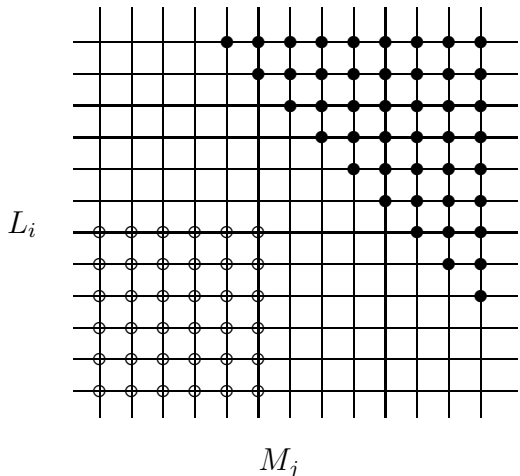
\begin{picture}(500,200)
\linethickness{0.15mm}
\multiput (130,30)(0,12){12}{\line(1,0){168}}
\multiput (140,20)(12,0){13}{\line(0,1){155}}
\multiput(133,27)(0,12){6}{ $\circ$}
\multiput(145,27)(0,12){6}{ $\circ$}
\multiput(157,27)(0,12){6}{ $\circ$}
\multiput(169,27)(0,12){6}{ $\circ$}
\multiput(181,27)(0,12){6}{ $\circ$}
\multiput(193,27)(0,12){6}{ $\circ$}
\multiput(281,159)(-12,0){9}{$\bullet$}
\multiput(281,147)(-12,0){8}{$\bullet$}
\multiput(281,135)(-12,0){7}{$\bullet$}
\multiput(281,123)(-12,0){6}{$\bullet$}
\multiput(281,111)(-12,0){5}{$\bullet$}
\multiput(281,99)(-12,0){4}{$\bullet$}
\multiput(281,87)(-12,0){3}{$\bullet$}
\multiput(281,75)(-12,0){2}{$\bullet$}
\multiput(281,63)(-12,0){1}{$\bullet$}
\put(105,90){$L_i$}
\put(200,0){$M_j$}
\end{picture}

There are $r$ solid dots along the top line, decreasing down to 1 on the $r$-th line; these dots  represent the lines that form the subcurve $Y$ required in (\ref{subcurve req}) above.  It is not hard, using methods of \cite{MN3} or \cite{MN2}, to check that $Y$ is arithmetically Cohen-Macaulay with maximal $h$-vector $(1,2,\dots,r)$.  The open dots represent the residual complete intersection of type $(t,t)$ (in our case, with $r \equiv 0 \hbox{ (mod 3)}$, we have $t = \frac{2r}{3}$).  Hence the dots that are not open represent our curve $C$, and clearly it has the required subcurve.

Now, observe that the initial degree of $I_Y$ is $r$, and that it contains $r+1$ minimal generators of degree $r$ (this can be read immediately from the $h$-vector).  By looking at suitable products of $L_i$ and of $M_j$ we can link $Y$ to a residual curve $Y'$ inside a complete intersection of type $(r,r+1)$ (it does not matter if $Y'$ overlaps with the open dots).  Note that $Y'$ has the same $h$-vector as $Y$, and that the complete intersection is a stick figure.  Hence the sum $I_Y + I_{Y'}$ is the saturated ideal of a reduced, arithmetically Gorenstein set of points $X_1$ (as was done in \cite{MN3}) with $h$-vector
\[
\left ( 1,3,6,\dots, \binom{r}{2}, \binom{r+1}{2} ,\binom{r}{2}, \dots, 6,3,1 \right )
\]
that lies on $C$ (in particular).

To compute the $h$-vector of $X_2$ (the quadric hypersurface section of $C$), we ``integrate" the $h$-vector $\Delta^2 h_C$ of $C$ and then take suitable differences.    We obtain (focusing on the important part of the computation)

 {\footnotesize 
\[
\begin{array}{l|ccccccccccccccccccccccccccccccccccc}
\deg & 0 & 1 & 2 & 3  &  \dots & \frac{4r}{3}-2 & \frac{4r}{3}-1 & \frac{4r}{3} &  \dots  & 2r-3 &  2r-2 & 2r-1 &
2r  \\
\hline
\Delta^2 h_{C} & 1 & 2 & 3 & 4 & \dots & \frac{4r}{3}-1 & \frac{4r}{3} & \frac{4r}{3} &\dots   & 6 & 4 & 2\\
\Delta h_{C} & 1 & 3 & 6 & 10 &  \dots & \binom{4r/3}{2}  & \binom{(4r/3)+1}{2} & \binom{(4r/3)+2}{2}-1  & \dots  & e-6 & e-2 & e & e  \\ 
\Delta h_{X_2}& 1 & 3 & 5 & 7 & \dots & \frac{8r}{3}-3 & \frac{8r}{3}-1 & \frac{8r}{3} &  \dots  &&10 & 6 & 2
\end{array}
\]
}
\noindent (where $e$ is the degree of $C$, whose precise value is not important here).  
Note that we are assuming that $d \geq 12$, so $r \geq 6$ and $r+1 \leq \frac{4r}{3}-1$.
We now have all the ingredients, and we again use basic double G-linkage, with the modifications mentioned already.  We obtain the computation

 {\scriptsize 
\[
\begin{array}{l|ccccccccccccccccccccccccccccccccccc}
\deg & 0 & 1 & 2 & 3 & 4 &  \dots & r & r+1 & r+2 &  \dots   & 2r-2 & 2r-1 &
2r  \\
\hline
\Delta h_{X_1} & & & 1 & 3 & 6 & \dots & \binom{r}{2} & \binom{r+1}{2} & \binom{r}{2}&  \dots & 6 & 3 & 1 \\
\Delta h_{X_2} & 1 & 3 & 5 & 7 & 9 & \dots & 2r+1 & 2r+3 &  2r+5-\delta & \dots & 10 & 6 & 2  \\ \hline
\Delta h_Z & 1 & 3 & 6 & 10 & 15 & \dots & \binom{r+2}{2} & \binom{r+3}{2} & \binom{r+2}{2}+4-\delta & \dots & 16 & 9 & 3
\end{array}
\]
}

\noindent where $\delta = 1$ if $r = 6$ and $\delta = 0$ otherwise.  The important thing to notice is that from degree $r+1$ to $r+2$, $\Delta h_{X_1}$ is decreasing, $\Delta h_{X_2}$ is increasing, and $\Delta h_Z$ is decreasing.  Hence because of the latter we need surjectivity if WLP is to hold.  To show that we in fact do not have surjectivity, by  Proposition \ref{WLP criterion} we have to verify that a sufficiently general line $\lambda$ does not impose $r+3$ independent conditions on $(I_Z)_{r+2}$.  (By ``sufficiently general" we mean that the first linear form only has to miss the points of $Z$, and the second is general.)  We focus on the quadric cutting out $X_2$, and more precisely on the two points of intersection of $\lambda$ with this quadric.  It is not hard to see that when $r=6$, one of the two points imposes a new condition but the second does not; when $r>6$ neither point imposes a new condition. 

The other cases, $r \equiv 1 \hbox{ (mod 3)}$ and $r \equiv 2 \hbox{ (mod 3})$ are similar and are left to the reader.
\end{proof}

\begin{remark}
The reader will note from the last example that we can have consecutive values of the $h$-vector of $Z$ such that the second is strictly smaller than the first, and still have multiplication by a general linear form fail to be surjective.
\end{remark}

\newpage

%%%%%%%%%%%%%%%%%%%%%%%%%%%%%%%%%%%%%%%%%%%%%%%%

\section{Non-unimodality}

In the paper \cite{zanello}, Zanello proved the following theorem:

\begin{theorem}[\cite{zanello}, Theorem 3]
There exist non-unimodal level $h$-vectors of codimension 3.
\end{theorem}

\noindent He did this by constructing an explicit example (Example 2 of his paper).   In this section we show how, despite the very algebraic nature of Zanello's non-unimodal example, we can mimic his essential numerical data geometrically to produce a level set of points with non-unimodal $h$-vector.  The $h$-vector of our example is exactly the same as the Hilbert function of Zanello's example.

\begin{example} \label{zanello nonunimodality}
Let $C \subset \PP^3$ be an arithmetically Cohen-Macaulay curve with $h$-vector
\[
(1,2,3,4,5,6,7,8,9,10).
\]
That is, $C$ has degree 55 and arithmetic genus 276.  We further assume that $C$ contains as a subset an arithmetically Gorenstein set of points, $X_1$, with $h$-vector 
\[
(1,3,6,10,6,3,1).
\]
This can be constructed, for example, as follows.  We start with an arithmetically Cohen-Macaulay curve $C_1$ with $h$-vector $(1,2,3,4)$, and link to a similar curve $C_2$ via a complete intersection of type $(4,5)$.  The ideal $I_G := I_{C_1} + I_{C_2}$ is the saturated ideal of an arithmetically Gorenstein zero-dimensional scheme $X_1$ with the stated $h$-vector, which of course then lies on both $C_1$ and $C_2$.  By making general choices, we can arrange that $X_1$, $C_1$ and $C_2$ be reduced \cite{GM5}.  By using Liaison Addition we can add plane curves to $C_1$ to obtain an arithmetically Cohen-Macaulay curve $C$ with the stated $h$-vector, and since $C \supset C_1 \supset X_1$ we have the desired inclusion.  We do not need that $C$ be irreducible, but we do need that it be reduced.

Now we let $X_2$ be a general cubic hypersurface section of $C$. The $h$-vector of $X_2$ can be computed from that of $C$, and is given below.  Let $Z_1 = X_1 \cup X_2$.  We make the computation
 {\small 
\[
\begin{array}{l|ccccccccccccccccccccccccccccccccccc}
\deg & 0 & 1 & 2 & 3 & 4 & 5 & 6 & 7 & 8 & 9 & 10 & 11  \\
\hline
\Delta h_{X_1} & &&& 1 & 3 & 6 & 10 & 6 & 3 & 1  \\
\Delta h_{X_2} & 1 & 3 & 6 & 9 & 12 & 15 & 18 & 21 & 24 & 27 & 19 & 10 \\ \cline{1-13}
\Delta h_{Z_1} & 1 & 3 & 6 & 10 & 15 & 21 & 28 & 27 & 27 & 28 & 19 & 10 
\end{array}
\]
}

Now, both $X_1$ and $X_2$ are level, but $Z_1$ is not.  Its Artinian reduction has one-dimensional socle in degree 9 and 10-dimensional socle in degree 11.  However, we can find a subset, $Z$, of $Z_1$ with any truncated Hilbert function \cite{GMR}, and in particular, a subset with $h$-vector
\[
(1,3,6,10,15,21,28,27,27,28).
\]
Since $I_Z$ agrees with $I_{Z_1}$ in all degrees $\leq 8$, the only socle elements for the Artinian reduction of $R/I_Z$ occur in degree 9, i.e.\ $Z$ is level.  \qed
\end{example}

Hence we have shown the following:

\begin{theorem}
There exists a reduced, level set of points in $\PP^3$ whose Artinian reduction has non-unimodal  Hilbert function.
\end{theorem}

\begin{remark}  Zanello also gives a more general result, namely that Artinian level algebras can be found which not only fail to be unimodal, but in fact have arbitrarily many ``valleys."  We have tried to reproduce this result in the context of points in $\PP^3$.  Unfortunately, our methods have thus far been unable to extend the result beyond one ``valley."
\end{remark}

\begin{remark}
Zanello mentions that Iarrobino has found an example, and verified on the computer, of an Artinian, codimension 3 level algebra with socle type 6 and non-unimodal Hilbert function.  We have not been able to find any example other than the one given above, in the context of reduced sets of points in $\PP^3$.  We have not seen Iarrobino's example, but perhaps seeing his approach would suggest a way to apply our methods to find reduced sets of points with non-unimodal $h$-vector and smaller socle type.

\end{remark}

%%%%%%%%%%%%%%%%%%%%%%%%%%%%%%%%%%%%%%%%%%%%%%%%%

\section{Other Unusual Behavior} \label{other unusual beh}

In this section we exhibit other unusual behavior that can be obtained for a level reduced set of points, answering some natural questions.  We also show some limitations of our method.

We begin by showing that even for reduced sets of points, it is not true that WLP implies SLP.

\begin{example}
Let $C_1 \subset C_2$ be arithmetically Cohen-Macaulay curves with $h$-vectors
\[
(1,2,3,3,2,1) \ \ \ \  \hbox{and} \ \ \ \ (1,2,3,4,5,5,5),
\]
respectively.  $C_1$ will be a complete intersection of type $(3,4)$.  To produce these curves, we will again use our trick of beginning with points in $\PP^2$ and taking a cone.  So let $Z_1$ be a general complete intersection of type $(3,4)$ in $\PP^2$.  Let $X$ be a general set of eight points in $\PP^2$.  Then $Z_1 \cup X$ has $h$-vector
\[
(1,2,3,4,5,5).
\]
Let $F$ be a plane curve of degree $5$ containing $Z_1 \cup X$.  Let $Y$ be the complete intersection set of points cut out by $F$ and a general quadric.  Then $Z_2 := (Z_1 \cup X) \cup Y$ is a basic double link, and has the desired $h$-vector given above.  Choosing a general point in $\PP^3$ and taking cones as before, we obtain the desired curves $C_1, C_2$.  

Let $X_1, X_2$ be general hyperplane sections (with different hyperplanes!) of $C_1, C_2$ and let $Z := X_1 \cup X_2$.  Then as before, $Z$ is obtained by basic double G-linkage, and its $h$-vector is computed by
 {\small 
\[
\begin{array}{l|ccccccccccccccccccccccccccccccccccc}
\deg & 0 & 1 & 2 & 3 & 4 & 5 & 6   \\
\hline
\Delta h_{X_1} && 1 & 2 & 3 & 3 & 2 & 1   \\
\Delta h_{X_2} & 1 & 2 & 3 & 4 & 5 & 5 & 5 \\ \cline{1-8}
\Delta h_Z & 1 & 3 & 5 & 7 & 8 & 7 & 6 
\end{array}
\]
}

We claim (and have verified on CoCoA) that $Z$ has WLP, but the multiplication from degree 3 to degree 5 by a general quadric is not surjective, and hence $Z$ does not have SLP.  We have to check that  (1) WLP holds (rather than fails, as has been the case until now in this paper), and (2)  SLP fails.     

For (1) we defer to CoCoA, since our methods for checking injectivity are not as effective, but it could also be checked by hand.  Once CoCoA verifies that a specific linear form gives the predicted multiplication, then it is true for the general linear form.   For (2) it is not enough to check that a ``random" form produced by CoCoA fails to give the expected multiplication.  Instead,  we take $\lambda$ to be a general complete intersection of type $(1,2)$ rather than $(1,1)$, reflecting that first we have to reduce by a linear form to get the Artinian reduction, and then we have to multiply by a general quadric.  Note that  again it is a question of the failure to impose independent conditions.  To check that SLP does not hold, we note that the first point of intersection of the conic $\lambda$ with the plane of $X_2$ does impose a new condition, but then the second point does not.  So we do not get the expected number of conditions and SLP fails, exactly as WLP failed in our other situations.
We leave the details to the reader, although we point out that the computation above suggests that multiplication by a general linear form should be of maximal rank (by imagining a decomposition of the form suggested by the diagram), but multiplication by a general quadric should have rank at most $2+4=6$, rather than 7, from degree 3 to degree 5.  \qed
\end{example}

As a result of this example, we have shown:

\begin{proposition}
There exist reduced, level zero-dimensional schemes whose general Artinian reduction has WLP but not SLP.
\end{proposition}

The following result shows that the value of the $h$-vector of our constructed set of points can be constant for arbitrarily many consecutive values (but the value gets correspondingly bigger), and WLP can correspondingly fail for arbitrarily many consecutive degrees.

\begin{proposition} \label{any d}
Let $d \geq 2$ be a positive integer.  Then there exists a reduced, level set of points, $Z \subset \PP^3$, whose $h$-vector has a constant value in each of $d$ consecutive degrees, and for which the corresponding $d-1$ maps  in {\em any} Artinian reduction all fail to be surjective (and hence WLP fails in each of these degrees).
\end{proposition}

\begin{proof}
We choose a reduced complete intersection $Y_1 \subset \PP^2$ of type $(d,d)$, and extend it to a set of points $Y_2 \subset \PP^2$ by adding $\binom{2d+1}{2} - d^2 = d(d+1)$ {\em general} points.  Choose a general point $P \in \PP^3$ and let $C_1, C_2$ be the cones over $Y_1,Y_2$ with vertex $P$.    Let $X_1,X_2$ be general hyperplane sections (different hyperplanes) of $C_1$ and $C_2$, respectively, and let $Z = X_1 \cup X_2$.  The $h$-vector of $Z$ is computed by 

 {\small 
\[
\begin{array}{l|ccccccccccccccccccccccccccccccccccc}
\deg & 0 & 1 & 2 & 3 & \dots &  d-1 & d & d+1 & \dots & 2d-3  & 2d-2 & 2d-1 \\
\hline
\Delta h_{X_1} && 1 & 2 & 3 & \dots & d-1 & d &  d-1 & \dots & 3 & 2 & 1 \\ 
\Delta h_{X_2} & 1 & 2 & 3 & 4 & \dots & d & d+1 & d+2 & \dots &  2d-2 & 2d-1 & 2d   \\ \hline
\Delta h_Z & 1 & 3 & 5 & 7 & \dots & 2d-1 & 2d+1 & 2d+1 & \dots & 2d+1 & 2d+1 & 2d+1
\end{array}
\]
}

\noindent We see that there are $d$ consecutive degrees with value $2d+1$.  The fact that each of the corresponding $d-1$ maps (multiplication by a general linear form)  in the Artinian reduction all fail to be surjective comes as before, using Proposition \ref{WLP criterion} and the fact that the point of intersection of $\lambda$ with the plane of $X_2$ is forced to be part of the base locus in the corresponding degree, and hence fails to impose an additional condition.
\end{proof}

\begin{remark}
The first case of Proposition \ref{any d}, namely the case $d=2$, which produces the $h$-vector $(1,3,5,5)$, should be isomorphic to the distraction (lifting) of Zanello's example coming from a monomial ideal.   We have not verified this, but numerically the examples are the same.  See also Example \ref{Zanello's points}.
\end{remark}

\begin{example} \label{limitations}
One of the cases left open in this paper is whether there exists a level reduced set of points whose Artinian reduction has type 2 but does not have WLP.  We now give a construction that, according to the ``philosophy" described in the introduction, {\em should} fail to have WLP, but does not.  Consider a reduced arithmetically Cohen-Macaulay curve $C_1$ in $\PP^3$ with $h$-vector $(1,2,3)$; it has degree 6 and arithmetic genus 3.  $C_1$ can be linked by a sufficiently general complete intersection $C$ of type $(3,4)$ to another curve, $C_2$, of the same degree and genus.  The intersection of $C_1$ and $C_2$ is a reduced, arithmetically Gorenstein set of points $X_1$ \cite{GM5} with $h$-vector $(1,3,6,3,1)$.  On the other hand, $X_1$ lies on the complete intersection curve $C$.  Taking a cubic hypersurface section of $C$ we obtain a complete intersection set of points, $X_2$.  The union $Z = X_1 \cup X_2$ is level of type 2, and has Hilbert function computed by

 {\small 
\[
\begin{array}{l|ccccccccccccccccccccccccccccccccccc}
\deg & 0 & 1 & 2 & 3 & 4 & 5 & 6 & 7  \\
\hline
\Delta h_{X_1} &&&& 1 & 3 & 6 & 3 & 1    \\
\Delta h_{X_2} & 1 & 3 & 6 & 8 & 8 & 6 & 3 & 1 \\ \cline{1-9}
\Delta h_Z & 1 & 3 & 6 & 9 & 11 & 12 & 6 & 2
\end{array}
\]
}

In order for the Artinian reduction of this set of points to have WLP, in particular the homomorphism from degree 4 to degree 5 induced by multiplication by a general linear form must be injective, i.e.\ have rank 11.  From the successful applications of this method in the previous sections, we would hope that the contribution to the rank from $X_1$ would be 3 and from $X_2$ to be 6, combining for a rank of 9 rather than 11.  However, CoCoA confirms that indeed the rank is 11, and in fact that the Artinian reduction of $Z$ does have WLP.  The problem comes from the fact that the three points of intersection of the general line $\lambda$ have to fail to impose independent conditions on forms of degree 7, and this degree is too large for the geometric type of obstructions that helped us before.  We have found many such examples.  \qed
\end{example}

\begin{example} \label{not injective}
We now illustrate that Proposition \ref{injectivity criterion} can also be used in a geometric way, this time  to prove the failure of {\em injectivity} for a level algebra.  However, we have not found as many useful applications as we did for Proposition \ref{WLP criterion}.

Let $C_1 \subset \PP^3$ be a reduced complete intersection curve of type $(2,2)$, and let $C_2$ be a reduced complete intersection curve of type $(5,5)$ containing $C_1$.  Let $X_1$ be a general hyperplane section of $C_1$, cut out by a general linear form $L$, and let $X_2$ be a general quadric section of $C_2$, cut out by a general quadric $Q$.  Let $Z' = X_1 \cup X_2$.  The $h$-vector of $Z'$ is computed by

 {\small 
\[
\begin{array}{l|ccccccccccccccccccccccccccccccccccc}
\deg & 0 & 1 & 2 & 3 & 4 & 5 & 6 & 7 & 8 & 9  \\
\hline
\Delta h_{X_1} &&& 1 & 2 & 1    \\
\Delta h_{X_2} & 1 & 3 & 5 & 7 & 9 & 9 & 7 & 5 & 3 & 1 \\ \hline
\Delta h_{Z'} & 1 & 3 & 6 & 9 & 10 & 9 & 7 & 5 & 3 & 1
\end{array}
\]
}

Using the machinery described above, we note first that $Z'$ is not level: it has socle in degree 4 and in degree 9.  However, as before (Example \ref{zanello nonunimodality}), we can choose a subset, $Z$, of $Z'$ by truncation, with $h$-vector $(1,3,6,9,10)$, and this will be level of type 10.  We will show that for any Artinian reduction of $R/I_Z$,  the multiplication by a general linear form from degree 3 to degree 4 is not injective, thereby showing that WLP does not hold.  Since $R/I_Z$ agrees with $R/I_{Z'}$ in degrees $\leq 4$, it is enough to show this for $R/I_{Z'}$.  

By Proposition \ref{injectivity criterion}, we have to show that for a sufficiently general line $\lambda$, $(I_{Z'} \cdot I_\lambda)_4 \neq (I_{Z'} \cap I_\lambda)_4$ (where $I_\lambda = (L_1,L_2)$ and  $L_1$ is chosen simply to avoid the points of $Z'$, while $L_2$ is chosen generically).  Since $\dim (I_{Z'})_3 = 1$ (with basis $LQ$), we note that $(I_{Z'} \cdot I_\lambda)_4$ is generated by $(LQ \cdot L_1, LQ \cdot L_2)$.  So it is enough to show that $\dim (I_{Z'} \cap I_ \lambda)_4 > 2$.  By considering multiplication by $Q$, we see that it is enough to show that $\lambda \cup X_1$ lies on at least three independent quadrics.  But $\lambda$ alone lies on 7 quadrics, and $X_1$ imposes at most 4 additional conditions, so we are finished.  \qed
\end{example}

%%%%%%%%%%%%%%%%%%%%%%%%%%%%%%%%%%%%%%%%%%%%%%%%

\section{Open questions}

We remark that our methods become more and more difficult to apply as the socle type decreases.  In this section we list some questions that remain open, having to do with low socle type.

In his paper \cite{zanello}, Zanello asks (Question 6) ``What is the maximum type $t_0$ such that all the codimension 3 level $h$-vectors of type $t \leq t_0$ are unimodal?  In particular, is there always unimodality for $t=2$?"  The following questions may be viewed as a continuation of Zanello's line of questioning.

In codimension 3, it is not known whether every Artinian Gorenstein algebra has the WLP, although it is known that they all have unimodal Hilbert function \cite{stanley}.  The only result to date in this direction comes from \cite{HMNW}, where it is shown that every height three complete intersection has WLP.  It is also not known if there is an Artinian level algebra of type 2 that fails to have WLP.  Of course as a result it is not known if there is a reduced set of points whose Artinian reduction has these properties, although it is known \cite{GM5} that for any set of graded Betti numbers that occurs for an Artinian Gorenstein algebra of codimension 3, there is a reduced set of points in $\PP^3$ with these graded Betti numbers.  One can also ask, for any $t_0$, what is the smallest socle degree for which there is an Artinian level algebra (resp.\ a level set of points) with socle type $t_0$ and which fails to have WLP or fails to have unimodal $h$-vector.

In codimension 4 it is was shown by Ikeda \cite{ikeda} that not all graded Artinian Gorenstein algebras have WLP.  It is not known if they all have unimodal Hilbert function.  It is also not known if Ikeda's example lifts to points, or more generally if there is a reduced set of points in $\PP^4$ whose Artinian reduction fails to have WLP, even if the Hilbert function is unimodal.

In codimension 5 and higher it is known that not all graded Artinian Gorenstein algebras have unimodal Hilbert function (see e.g.\  \cite{ber-iar}, \cite{boij}, \cite{BL}).  But again, it is not known if this extends to points.  It would be extremely interesting to know if all arithmetically Gorenstein, reduced sets of points have (at least) the general Artinian reduction with WLP, because then thanks to \cite{MN3} we would have a complete description of the Hilbert functions of reduced, arithmetically Gorenstein sets of points.

We do not know of any example of a level set of reduced points where the general Artinian reduction has WLP but a special one does not.

It would be very interesting to know if we can use these (or related) methods to show that the postulation Hilbert scheme is reducible in some cases, for level sets of points.  Basic double linkage was used in \cite{families} to show the reducibility of the postulation Hilbert scheme in a different setting -- there a fundamental role was in fact played by  {\em forcing} the existence of socle elements in degrees other than the last one, so those results do not apply at all to our situation.

%%%%%%%%%%%%%%%%%%%%%%%%%%%%%%%%%%%%%%%%%%%%%%%%%%%%%%%%%%%%%%%%%%%%%%%%%%%
%%%%%%%%%%%%%%%%%%%%%%%%%%%%%%%%%%%%%%%%%%%%%%%%%%%%%%%%%%%%%%%%%%%%%%%%%%%

\end{document}